\documentclass{amsart}       
\usepackage{graphicx}
\usepackage{pgf, tikz,amsfonts,amsmath}
\usetikzlibrary{trees}
\newcommand{\sN}{\mathcal N}
\newcommand{\sY}{\mathcal Y}

\newcommand{\N}{\ensuremath{\mathbb N}}
\newcommand{\Ew}{\ensuremath{\mathcal P}}
\newcommand{\Ww}{\ensuremath{\mathcal W}}
\newcommand{\R}{\ensuremath{\mathbb{R}}}
\newcommand{\ew}{\ensuremath{\mathbf p}}
\newcommand{\ww}{\ensuremath{\mathbf w}}

\newcommand{\E}[1]{\ensuremath{\mathbf{E} \left[#1 \right]}}
\newcommand{\Prob}[1]{\ensuremath{\mathbb{P} \left(#1 \right)}}
\newcommand{\I}[1]{\ensuremath{\mathbf{1}_{  #1  }}}
\newcommand{\Law}{\mathcal{L}}

\newtheorem{theorem}{Theorem}

\newtheorem{lemma}{Lemma}
\newtheorem{proposition}{Proposition}

\author{Rafik Aguech}
 \address{Rafik Aguech, Department of Statistic and Operation Research, College of Sciences, King Saud University, P.O. Box 2455, 
Riyadh 11451, Saudi Arabia} \email{rafik.aguech@ksu.edu.sa}                                
           \author{Anis Amri} 
            \address{Anis Amri, University of Monastir, Avenue Taher Hadded B.P 56, Monastir 5000, Tunisia}  \email{anis.amri@hotmail.com}
          \author{Henning Sulzbach}
          \address{Henning Sulzbach, McGill University, 3480 University Street, H3A 0E9 Montreal, QC, Canada} \email{henning.sulzbach@gmail.com}
          \curraddr{School of Mathematics, University of Birmingham, Birmingham B15 2TT, Great Britain}

\title{On weighted depths  in random binary search trees}


\begin{document}
\begin{abstract}
Following the model introduced by Aguech, Lasmar and Mahmoud [\emph{Probab. Engrg. Inform. Sci.} \textbf{21} (2007) 133--141], the weighted depth of a node in a labelled rooted tree is the sum of all labels on the path connecting the node to the root. We analyze weighted depths  of nodes with given labels, the last inserted node, nodes ordered as visited by the depth first search process, the weighted path length and the weighted Wiener index in a random binary search tree. We establish three regimes of nodes depending on whether the second order behaviour of their weighted depths follows from fluctuations of the keys on the path, the depth of the nodes, or both.  Finally, we investigate a random distribution function on the unit interval arising as scaling limit for weighted depths of nodes with at most one child.
\end{abstract}

\date{\today}
\maketitle

\section{Introduction}
\label{intro}
\begin{figure} [htb] \label{fig_ex}
\centering
\begin{tikzpicture}[level distance=1.5cm,
level 1/.style={sibling distance=3.5cm},
level 2/.style={sibling distance=2cm}]
\tikzstyle{every node}=[circle,draw]

\node (Root) {4}
    child {
    node {2}
    child { node {1} }
    child {node {3} }
  }
child {
    node {6}
    child { node {5} }
    child { node {7} }
};
\end{tikzpicture}
\caption{binary search tree constructed from the list $4, 2, 6, 5, 7, 3,  1$.}
\end{figure}
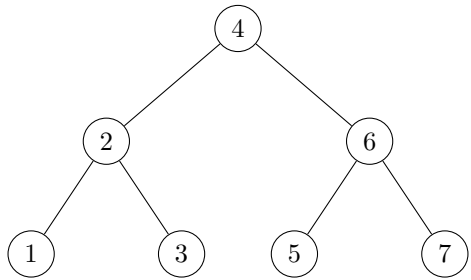

The binary search tree is an important data structure in computer science allowing for efficient execution of database operations such as insertion, deletion and retrieving of data. Given a list of elements $x_1, x_2, \ldots, x_n$ from a totally ordered set, 
it is the unique labelled rooted binary tree with $n$ nodes constructed by successive insertion of all elements  satisfying the following property: for each node in the tree with label (or key), say $y$, all keys stored in its left (right) subtree are at most equal to (strictly larger than) $y$.
For an illustration, see Figure \ref{fig_ex}.

Properties of binary search trees are typically analyzed under the random permutation model where the data $x_1, \ldots, x_n$ are generated by a uniformly chosen permutation of the first $n$ integers. Among the quantities studied in binary search trees, one finds depths of and distances between nodes related to the performance of search queries and finger searches in the database, the (total) path length measuring the cost of constructing the tree as well as the Wiener index. Further, more complex parameters such as  the height corresponding to worst case search times, the saturation level and the profile have been studied thoroughly. We review literature relevant in the context of our work below.

In this note we complement the wide literature on random binary search trees by the analysis of depths of nodes, path length and Wiener index in their weighted versions as introduced by Aguech, Lasmar and Mahmoud \cite{aglama07}. Here, the weighted depth of a node is the sum of all keys stored on the path to the root. In \cite{aglama07}, results about weighted depths of extremal paths have been obtained. Kuba and Panholzer \cite{kupa07}, \cite{kupa08} studied the problem in random increasing trees covering the random recursive tree and the random plane-oriented recursive tree.
Weighted depths of nodes and the weighted height were also studied by Broutin and Devroye \cite{brde} in a more general tree model, which relies on assigning weights to the edges of the tree.
Further, the weighted path length in this model was investigated by R{\"u}schendorf and Schopp \cite{ruschopp}.
Note that we deviate from the notation introduced in \cite{aglama07} and \cite{kupa07} using the term weighted depth for what is called weighted path length there since we also study a weighted version of the (total) path length of binary search trees.

\section{Preliminaries} \label{sec:not}

We introduce some notation. By the \emph{size} of a finite binary tree, we refer to its number of nodes. 
Upon embedding a finite rooted binary tree in the complete infinite binary tree, a node is called \emph{external} if its graph distance to the binary tree is one.
Any node on level $k \geq 1$ in a rooted binary tree is associated a vector $v_1 v_2 \ldots v_k \in \{0,1\}^k$ where $v_i = 0$ if and only if the path from the root to the node continues in the left subtree upon reaching level $i-1$. 

Let $n \geq 1$ and $1 \leq k \leq n$. Under the random permutation model (short: \emph{permutation model}), let $D_k(n)$ be the depth of the node labelled $k$. 
By $W_k(n)$ we denote the sum of all keys on the path from the root to the node labelled $k$ including the labels of both endpoints. 
For $x = x_1 x_2 \ldots \in \{0,1\}^\infty$,  let $B_n(x)$ be the maximal depth among nodes of the form $x_1 \ldots x_k, k \geq 0$.
We use $X_n$ ($\mathbb X_n$) to denote the (weighted) depth of the $n$th inserted node. Finally, we define the height of the tree by 
 $H_n = \sup \{k \in \N : D_k(n) > 0\}$.

Throughout the paper, we denote by $\Law(X)$ the distribution of a random variable $X$. For real-valued $X$ with finite second moment, we write $\sigma_X$ for its standard deviation. 
By $\sN$ we denote a random variable with the standard normal distribution, and by $\mu$ the Dickman distribution on $[0, \infty)$ characterized by its Fourier transform,
\begin{align} \label{fourier_dickman} \int e^{i \lambda x} d \mu(x) = \exp \left( \int_0^1 \frac{e^{i \lambda x}-1}{x} dx\right), \quad \lambda \in \R. \end{align}
The origins of the Dickman distribution go back to Dickman's \cite{dickman1930} classical result on large prime divisors. Compare Hildebrandt and Tenenbaum \cite{hildetenen} for a survey on the problem. 
In the probabilistic analysis of algorithm, $\mu$ first arose in Hwang and Tsai's \cite{hwangtsai} study of the complexity of Hoare's selection algorithm. 
We refer to this work for a discussion of more details on the distribution, historical background and further references.

Finally, we use the Landau notations little--$o$, big--$O$, little--$\omega$, big--$\Omega$ and big--$\Theta$ as $n \to \infty$. 

\subsection{Depths and height} \label{sec:dh} We recall the following fundamental property of random binary search trees going back to Devroye \cite{devroyeheight}: 
 in probability and with respect to all moments, we have
\begin{align} \label{conv_height} \frac{H_n}{\log n} \to c^*, \end{align}
where $c^* = 4.31 \ldots$ is the larger of the two solutions to the transcendent equation $e = (\frac{2e}{x} )^x$.  
Next, by classical results due to  Brown and Shubert \cite{brownshu} and Devroye \cite{dev88}, for any $x \in \{0,1\}^\infty$, in distribution,
\begin{align} \label{dev2} \frac{B_n(x) - \log n}{\sqrt{\log n}} \to \sN, \quad \frac{X_n - 2 \log n}{\sqrt{ 2\log n}} \to \sN. \end{align}
(In \cite[Theorem O1]{dev88}, the first convergence in the last display is formulated for $x = \mathbf{0} := 00\ldots$ The general case follows, since, by symmetry, $\mathcal L( B_n(x)) = \mathcal L( B_n (\mathbf 0)) $ 
for all $x$. The second convergence was also claimed in a footnote by Mahmoud and Pittel \cite{mapi84}.) 
Gr{\"u}bel \cite{gruebelsilhouette} studied the process $\{B_n(x): x \in \{0, 1 \}^\infty\} $, the so-called \emph{silhouette}, thereby obtaining a functional limit theorem for its integrated version.
The asymptotic behaviour of depths of nodes with given labels has been analyzed by Devroye and Neininger \cite{denefinger04}:
uniformly in $1\leq k \leq n$ and as $n \to \infty$,
\begin{align} \E{D_k(n)} = \log(k(n-k)) + O(1), \quad \text{Var}(D_k(n)) = \log(k(n-k)) + O(1). \label{expansionD} \end{align}
Moreover, for any $1 \leq k \leq n$, which may depend on $n$, in distribution
\begin{align}  \frac{D_k(n) - \E{D_k(n)}}{\sigma_{D_k(n)}} \to \sN.  \label{limDN} \end{align}
Here, one should also compare Gr{\"u}bel and Stefanoski \cite{gruebelpoisson05} for stronger results in the context of the corresponding Poisson approximation.
 For a survey on depths and distances in binary search trees, we refer to Mahmoud's book \cite{mahmoudbook}.
Finally, the asymptotic behaviour of the weighted depths of the nodes associated with the vectors $\mathbf 0$ and $\mathbf 1 := 11\ldots$ denoted by $\mathcal{L}_n$ and $\mathcal{R}_n$ ($\mathcal L$ and $\mathcal R$ 
stand for left and right) were studied  in \cite{aglama07}. In distribution, 
\begin{align} \label{lasm}
\frac{\mathcal{L}_n}{n} \to \sY, \quad \frac{\mathcal{R}_n - n B_n(\mathbf{1})}{n \sqrt{\log n}} \to 0, \end{align}
where $\sY$ has the Dickman distribution. The first convergence is closely related to the limit law in Theorem 3.1 in \cite{hwangtsai}.

\subsection{Path length and Wiener index} \label{sec:pw}
In a rooted tree, the path length is defined as the sum over all depths of nodes. Moreover, the Wiener index is obtained by summing  all distances of unordered pairs of vertices. For a random binary search tree of size $n$, we denote its path length by $P_n$ and its Wiener index by $W_n$.
Denoting by $\gamma$ the Euler-Mascheroni constant, 
we have
\begin{align} \label{meanpath}
 \E{P_n} = 2n \log n + (2\gamma -4 )n + o(n),  \quad \text{Var}(P_n) = \frac{21 - 2 \pi^2}{3} n^2 + o(n^2),  \end{align}
going back to Hoare \cite{hoare} and Knuth \cite{knuth1973b}. Further, by \cite{newiener}, \begin{align} \label{meanwiener} \E{W_n} = 2n^2 \log n + (2\gamma -6)n^2 + o(n^2), \quad  \text{Var}(W_n) = \frac{20 - 2 \pi^2}{3}  n^4 + o(n^4). \end{align}
Central limit theorems for the path length go back to
R{\'e}gnier \cite{regnier} and R{\"o}sler \cite{ro91}, for the Wiener index to Neininger \cite{newiener}. More precisely, by \cite[Theorem 1.1]{newiener}, there exists a non-trivial random variable $Z^*$ on $\R^2$ characterized by a stochastic fixed-point equation, such that, in distribution, 
 \begin{align}
\left(\frac{W_n-\mathbb{E}[W_n]}{n^2},\frac{P_n-\mathbb{E}[P_n]}{n}\right) \to Z^*. \label{lim_oldnei}
\end{align}

\subsection{The i.i.d.\ model} \label{sec:iid}
We also consider binary search trees of size $n$ where the data are chosen as the first $n$ values of a sequence of independent random variables $U_1, U_2, \ldots$ each having the uniform distribution on $[0,1]$. 
Since the vector $(\text{rank}(U_1), \ldots, \text{rank}(U_n))$ constitutes a uniformly chosen permutation, in distribution,  both the permutation model and the i.i.d.\ model lead to the same unlabelled tree. We use the same notation as in the permutation model 
for quantities not involving the labels of nodes, that is, $X_n, h_n, H_n, P_n, W_n$ and $B_n(x)$. Further, we define the weighted path length $\mathcal{P}_n$ as the sum of all weighted depths, and the weighted Wiener index $\mathcal{W}_n$ as the sum over all pairs of weighted distances. Here, the weighted distance between two nodes equals the sum of all labels on the path connecting them, labels of endpoints included. (Notice that the weighted distance between a node and itself is equal to its label.) Finally, analogously to $B_n(x)$, we define $\mathcal B_n(x)$ as the weighted depth of the node of largest depth on the path $x$. We call $\{\mathcal B_n(x): x \in \{0,1 \}^\infty\}$, the weighted silhouette of the tree (at time $n$).

\section{Main results} \label{sec main}
Our main results are divided into two groups: Theorems \ref{thm_large} and \ref{thm_small} hold in the permutation model while Theorems \ref{thm_binary} and \ref{thm_pathlength} are formulated in the i.i.d.\ model.

\subsection{Results in the permutation model}
We start with the expansions of the first two moments of the weighted depth $W_k(n)$. Uniformly in $1 \leq k \leq n$, as $n \to \infty$, 
\begin{align} 
& \E{W_k(n)} = k \log(k(n-k+1)) + n + O(k + \log n), \label{expansionW} \\ 
& \text{Var}(W_k(n)) = k^2 \log(k(n-k+1)) + \frac{n^2}{2} + O(kn). \label{expansionV} 
\end{align}
It turns out that the asymptotic distributional behavior of $W_k(n)$ with respect to terms of second order is entirely described by
that of $k D_k(n)$ if and only if $k = \omega(n / \sqrt{\log n})$.  \sloppy  Accordingly, in the remainder of this paper, we call nodes with labels of order \sloppy $\omega(n / \sqrt{\log n})$ \emph{large} and of order $O(n / \sqrt{\log n})$ \emph{small}. 

\begin{theorem}[Weighted depths of large nodes] \label{thm_large}
For $k = \omega(n / \sqrt{\log n})$,
\begin{align} \label{diff_WD}
 \E{| W_k(n) - k D_k(n)|} = o(\sigma_{k D_k(n)}).
\end{align}
In particular, for $0 < \alpha < 1$ and  $|k/n - \alpha| = o((\log n)^{-1/2})$, in distribution,
\begin{align} \label{lim_DW_large}
\left(\frac{D_k(n) - 2 \log n}{\sqrt{2 \log n}}, \frac{W_{k}(n) - 2\alpha n \log n}{\alpha n \sqrt{2 \log n}}  \right) \to (\sN,\sN). \end{align}
For the last inserted node, in distribution,
\begin{align} \label{lim_last}
\left(\frac{X_n - 2 \log n}{\sqrt{2 \log n}}, \frac{\mathbb X_n}{2n \log n} \right) \to \left(\sN, \xi \right), \end{align}
where $\sN$ and $\xi$ are independent and $\xi$ is uniformly distributed on $[0,1]$.
\end{theorem}
The asymptotic behavior of weighted depths of small nodes is to be compared with the corresponding results in \cite{kupa07}. Here, another phase transition occurs when $k = o(n / \sqrt{\log n})$.
\begin{theorem}[Weighted depths of small nodes] \label{thm_small}
Let $k = O(n/ \sqrt{\log n})$. Then, in distribution,
\begin{align} \label{lim_small} \left( \frac{D_k(n) - \E{D_k(n)} }{\sigma_{D_k(n)}}, \frac{W_k(n) - k D_k(n)}{n} \right) \to (\sN, \sY), \end{align}
where $\sN$ and $\sY$ are independent and $\sY$ has the Dickman distribution.
Thus, if $k \sqrt{\log n}/n \to \beta \geq 0$, in distribution,
$$\left(\frac{D_k(n) - 2 \log n }{\sqrt{2 \log n}}, \frac{W_{k}(n) - \E{W_k(n)}}{n}  \right) \to (\sN, \sY + \sqrt{2} \beta \sN -1).$$
In particular, if $|k \sqrt{\log n}/n - \beta| = o((\log n)^{-1/2})$ with $\beta > 0$, then, in distribution, 
$$\left(\frac{D_k(n) -  2\log n}{\sqrt{2 \log n}}, \frac{W_{k}(n) - 2 \beta n \sqrt{\log n}}{n}  \right) \to (\sN, \sY + \sqrt{2}\beta \sN).$$
\end{theorem}

 \subsection{Results in the i.i.d.\ model}
\begin{figure} [htb]  \begin{center} \includegraphics[width = 12cm, height = 6cm] {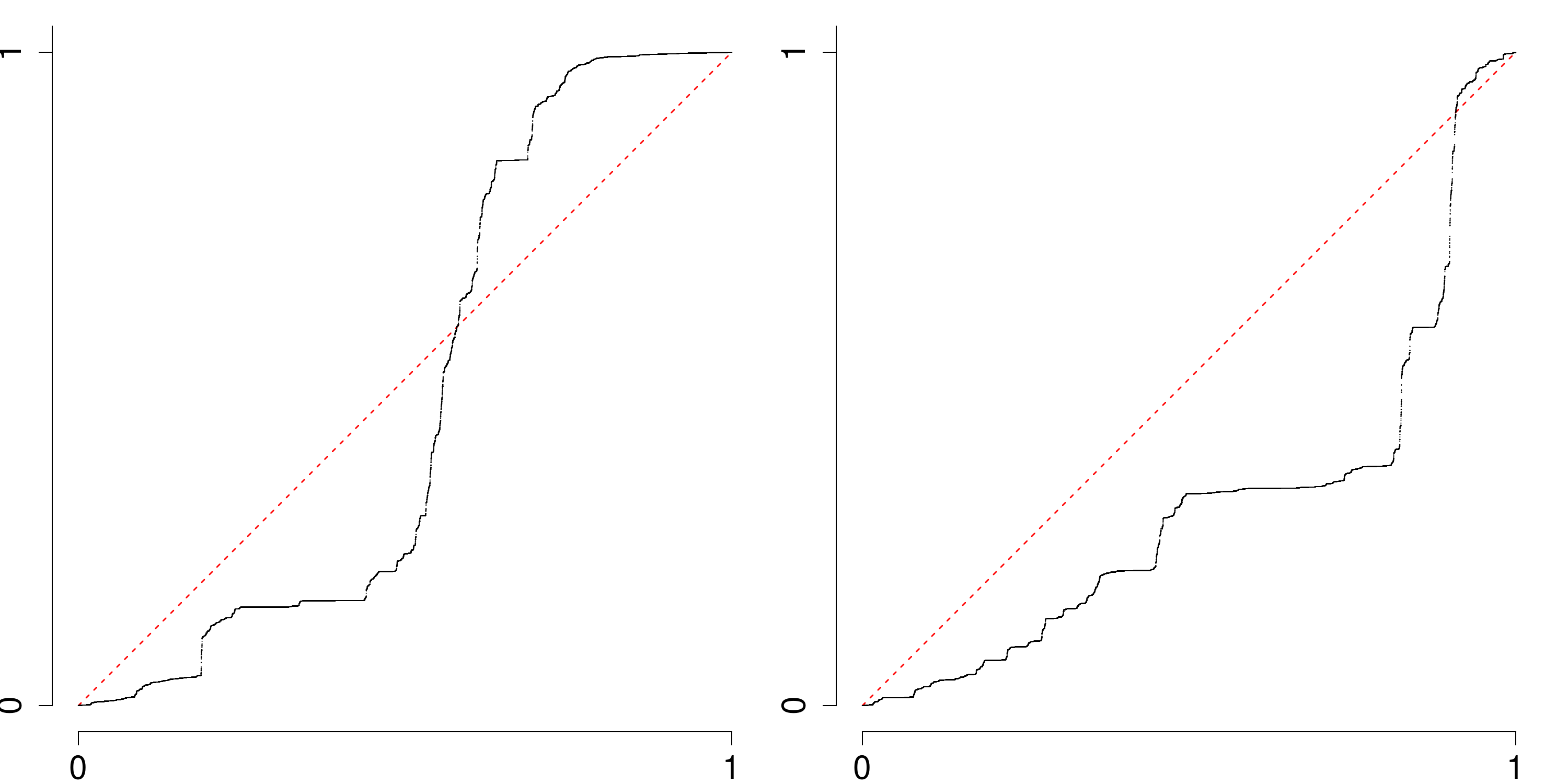}   \end{center}\caption{two simulations of $\Xi_{15}$, the dotted line being the graph of the identity function.}  \label{fig1} \end{figure}

Any $x \in \{0,1\}^\infty$ corresponds to a unique value $x \in [0,1]$ by $x = \sum_{i=0}^\infty x_i 2^{-i}$. This identification becomes one-to-one upon allowing only those $x \in \{0,1\}^\infty$ which contain infinitely many zeros and $x = \mathbf{1}$. In the i.i.d.\ model, 
for any $x \in \{0,1\}^{\infty}, k \geq 1$, the node $x_1 \ldots x_k$ eventually appears in the sequence of binary search trees and we write $\Xi_k(x)$ for its ultimate label. The following theorem about the behavior of $\mathcal{B}_n(x)$ involves a random continuous distribution function arising as the almost sure limit of $\Xi_k(x), x \in [0,1],$ as $k \to \infty$. We believe that this process is of independent interest and state some of its properties in Proposition \ref{propxi} in Section \ref{sec13}. The simulations of $\Xi_{15}$ presented in Figure \ref{fig1} illustrate the scaling limit.

\begin{theorem}[Weighted silhouette] \label{thm_binary}  There exists a random continuous and strictly increasing bijection $\Xi(x), x \in [0,1]$, such that, almost surely, uniformly on the unit interval, 
 $\Xi_k(x) \to \Xi(x)$.  For any $x \in [0,1]$,  in probability,  
\begin{align}  \label{conv_rep} \frac{\mathcal{B}_n(x)}{\log n} \to \Xi(x). \end{align}
Also, for any $m \geq 1$, in probability
\begin{align} \label{Lm_conv} \int_0^1 \left | \frac{\mathcal{B}_n(x)}{\log n} - \Xi(x) \right|^m dx \to 0. \end{align}
Further, in probability, \begin{align} \label{weighted_height} \sup_{x \in [0,1]} \frac{\mathcal{B}_n(x)}{\log n} \to c^* = 4.31\ldots \end{align} with $c^*$ as in \eqref{conv_height}.
Finally, for any $x \in [0,1]$, in distribution,
\begin{align} \label{conv_binary_joint}
\left( \frac{B_n(x) - \log n}{\sqrt{\log n}}, \frac{\mathcal{B}_n(x)}{\log n} \right) \to (\sN, \Xi(x)),
\end{align}
where $\sN$ and $\Xi(x)$ are independent. \end{theorem}

The next theorem extends the distributional convergence result in Theorem 1.1 in \cite{newiener}, that is \eqref{lim_oldnei},  by central limit theorems for the weighted path length and the weighted Wiener index.
\begin{theorem}[Weighted path length and Wiener index] \label{thm_pathlength}
In the i.i.d.\ model, we have 
$$ \E{\Ew_n}  = n \log n  + (\gamma - 3/2) n + o(n), \quad \E{\Ww_n} =  n^2 \log n + (\gamma -11/4)n^2 + o(n^2), $$
and 
\begin{align*}
\emph{Var}(\Ew_n) =   \frac{65 - 6 \pi^2}{36} n^2 + o(n^2), \quad \emph{Var}(\Ww_n) = \frac{2413 - 240 \pi^2}{1440} n^4 + o(n^4).
\end{align*}
The leading constants in the expansions of the covariances between $P_n, W_n, \Ew_n$ and $\Ww_n$ are given in \eqref{covs}--\eqref{cove}. (The leading constant for $\emph{Cov}(P_n, W_n)$ was already given in \cite{newiener}.)
As $n \to \infty$, with convergence in distribution and with respect to the first two moments in $\R^4$, we have
\begin{align*}
\left(\frac{\Ww_n-\mathbb{E}[\Ww_n]}{n^2},\frac{W_n-\mathbb{E}[W_n]}{n^2},\frac{\Ew_n-\mathbb{E}[\Ew_n]}{n},\frac{P_n-\mathbb{E}[P_n]}{n}\right) \to Z, 
\end{align*}
where the limiting distribution $\Law(Z)$ is the unique fixed-point of the map $T$ in \eqref{fixT}.
\end{theorem}

\textbf{Conclusions.} We have seen that there exist three types of nodes showing significantly different behavior with respect to their weighted depths.
By Theorem \ref{thm_large},  for $k = \omega(n/\sqrt{\log n})$,  second order fluctuations of weighted depths
are due to variations of the depth of nodes. 
In the second regime, when $k = \Theta(n/\sqrt{\log n})$, variations of weighted depths are determined by two independent contributions, one for the depths and one for the keys on the paths. Finally,  when $k = o(n/\sqrt{\log n})$ only fluctuations of labels on
paths influence second order terms of weighted depths. The third regime can be further subdivided with respect to the first order terms  of $W_k(n)$ and $k D_k(n)$: for $k = \omega(n / \log n)$, they coincide, for $k = \Theta(n/\log n)$, they are of the same magnitude, whereas, for $k = o(n/ \log n)$, they are of different scale.
By Theorem \ref{thm_binary}, the weighted silhouette behaves considerably different. Here, the lack of concentration around the mean leads to an interesting random distribution function on the unit interval as scaling limit.

\subsection{Further results and remarks} \label{sec13} \hfill\\

 \textbf{Model comparison.}   We decided to present Theorems \ref{thm_binary} and  \ref{thm_pathlength} in the i.i.d.\ model rather than in the permutation model since this allows for a stronger mode of convergence in 
\eqref{conv_rep}, \eqref{Lm_conv} and a clearer presentation of the proof of Theorem \ref{thm_pathlength}. 
In the i.i.d.\ model, denoting by $\mathcal{W}_{(k)}(n)$ the weighted depth of the node of rank $k$ among the first $n$ inserted keys, Theorems \ref{thm_large} and \ref{thm_small} remain valid upon replacing $W_k(n)$ by $n \mathcal{W}_{(k)}(n)$. 
Similarly, Theorems \ref{thm_binary} and \ref{thm_pathlength} hold in the permutation model where weighted depths and the weighted path length are to be scaled down by a factor $n$ and the weighted Wiener index by a factor $n^2$. 
The convergences in \eqref{conv_rep} and \eqref{Lm_conv} then only hold in distribution.
This can be deduced most easily from the following coupling of the two models: starting with the binary search tree in the i.i.d.\ model, also consider the random binary search tree in the permutation model relying on the permutation $(\text{rank}(U_1), \ldots, \text{rank}(U_n))$. Then, for all $1 \leq k \leq n$,
\begin{align} \left | \mathcal{W}_{(k)}(n) - \frac{W_k(n)}{n} \right| \leq H_n \max_{1 \leq i \leq n} \left| U_i - \frac{\text{rank}(U_i)}{n} \right|. \label{conn_models} \end{align}
It is well-known that the second factor on right hand side grows like $n^{-1/2}$, compare, e.g.\ Donsker's theorem for empirical distribution functions or the Dvoretzky-Kiefer-Wolfowitz inequality \cite{dvkiwo}. Combining this, \eqref{conn_models} and \eqref{conv_height} is sufficient to transfer all results in Section \ref{sec main} between the two models.

\medskip \textbf{The depth first search process.}
In the permutation model, let $v_1, \ldots, v_{n+1}$ be the external nodes as discovered by the depth first search process from left to right. By $D^*_k(n)$ and $W^*_k(n), 1 \leq k \leq n+1$, we denote depth and weighted depth of the external node $v_k$. Then, at the end of Section \ref{sec:typical}, we show that, uniformly in $1 \leq k \leq n$,
\begin{align} \label{Dbound}
\E{|D_k(n) - D^*_k(n)|^2} = o(\log n), \quad \E{|W_k(n) - W^*_k(n)|^2} = o(\text{Var}(W_k(n))).\end{align}
Thus, the results in Theorems \ref{thm_large} and \ref{thm_small} also cover the second order analysis of the sequences $D_k^*(n)$ and $W_k^*(n)$.

\medskip  \textbf{Weighted distances.}
In the permutation model, let $D_{k,\ell}(n)$ be the graph distance between the nodes labelled $1 \leq k \leq \ell \leq n$ and $W_{k,\ell}(n)$ be the sum of all labels on the path from $k$ to $\ell$, labels at the endpoints included. Asymptotic normality for the sequence $(D_{k,\ell}(n))$ (after rescaling) under the optimal condition $\ell - k \to \infty$ has been obtained in \cite{denefinger04}. For uniformly chosen nodes, distributional convergence results date back to Mahmoud and Neininger \cite{mane03} and Panholzer and Prodinger \cite{papr04}.
Analogously to Theorem \ref{thm_large}, it is straightforward to prove central limit theorems jointly for weighted and non-weighted distances. We only state the results. 
If $\ell - k = \Omega(n)$ and $k = \omega(n / \sqrt{\log n})$, then
\begin{align*} \E{|W_{k,\ell}(n) - kD_k(n) - \ell D_\ell(n)|} = \sigma_{\ell D_\ell(n)}. \end{align*}
In particular, for $0 < s < t < 1$ and $|k/n - s| = o((\log n)^{-1/2}), |\ell/n - t| = o((\log n)^{-1/2})$, we have, in distribution,
\begin{align*}
\Bigg(\frac{D_k(n) - 2 \log n}{\sqrt{2 \log n}}, \frac{D_\ell(n) - 2 \log n}{\sqrt{2 \log n}}, &  \frac{D_{k,\ell}(n) - 4 \log n}{\sqrt{4 \log n}},  \frac{W_{k,\ell}(n) - 2 (s+t) n \log n}{n \sqrt{2 \log n}}\Bigg) \nonumber \\
& \to \left(\sN_1, \sN_2, \frac{\sN_1 + \sN_2}{\sqrt{2}}, s \sN_1 + t \sN_2\right).   \end{align*}
Here, $\sN_1, \sN_2$ are  independent random variables both with the standard normal distribution. 

\medskip  \textbf{The limit process $\Xi$.} The process $\Xi$ in Theorem \ref{thm_binary} is a random distribution function. In particular, it can be regarded as an element in the set of c{\`a}dl{\`a}g functions $\mathcal{D}[0,1]$ consisting of all $f: [0,1] \to \R$, such that, for all $t \in [0,1]$, $f(t) = \lim_{s \downarrow t} f(s)$ and $\lim_{s \uparrow t} f(s)$ exists. The absolute value of $f$ is defined by $\sup \{|f(t)| : t \in [0,1]\}$. Endowed with Skorokhod's topology $J_1$, $\mathcal{D}[0,1]$ becomes a Polish space. We refer to Chapter 3 in  Billingsley's book \cite{bil99} for detailed information on this matter.
\begin{proposition}[Properties of $\Xi$] \label{propxi}
 The process $\Xi$ is unique (in distribution) among all c{\`a}dl{\`a}g processes with finite absolute second moment satisfying 
\begin{align} \label{fix_process}
\Law((\Xi(t))_{t \in [0,1]}) = \Law(\left( \I{[0,1/2)}(t) U \Xi(2t) +  \I{[1/2,1)}(t) \left( (1-U) \Xi'(2t-1) + U \right) \right)_{t \in [0,1]}).
\end{align} 
Here, $\Xi, \Xi', U$ are independent, $U$ has the uniform distribution on $[0,1]$ and $\Xi'$ is distributed like $\Xi$.
We have
\begin{enumerate} 
\item   $\E{\Xi(t)} = t$ for all $t \in (0,1)$;
\item  $\Law( (\Xi(t))_{t \in [0,1]} ) = \Law ( (1 - \Xi(1-t))_{t \in [0,1]})$;
\item $\Xi(\xi)$ has the arcsine distribution with density
$$\frac{1}{\pi \sqrt{x(1-x)}}, \quad x \in (0,1),$$
where $\Xi, \xi$ are independent and $\xi$ has the uniform distribution on $[0,1]$;
\item for $t \in (0,1)$, $\Law(\Xi(t))$ has a smooth density $f_t: (0,1) \to (0, \infty)$;
\item for $t \in (0,1/2)$,  $x f'_t(x) = - f_{2t}(x)$, $x \in (0,1)$, $f_t$ is strictly monotonically decreasing and $\lim_{x \uparrow 1} f_t(x) = 0$;
\item with $\alpha^{(i)}_t := \lim_{x \downarrow 0} f^{(i)}_t(x)$, $i = 0,1, t \in (0,1/2)$ and $\gamma_0 = 1/4, \gamma_1 = 5/16$,  we have $\alpha^{(i)}_t =  (-1)^i \infty$  for $0 < t \leq \gamma_i$, $|\alpha^{(i)}_t| < \infty$ for $\gamma_i < t < 1/2$ and $|\alpha^{(i)}_t| \uparrow  \infty$ as $t \downarrow \gamma_i$.
\end{enumerate}
\end{proposition}

\medskip  \textbf{Random recursive trees.}
A random recursive tree is constructed as follows: starting with the root labelled one, in the $k$th step, $k \geq 2$, a node labelled $k$ is inserted in the tree and connected to an already existing node chosen uniformly at random. Weighted depths in random binary search trees differ substantially from those in random recursive trees analyzed in \cite{kupa07} where all nodes show an asymptotic behaviour comparable to that  of nodes labelled $k = o(n/\sqrt{\log n})$ in the binary search tree. The difference is highlighted by the weighted path length. Being of the same order as the path length in binary search trees, it follows from results in \cite{kupa07} that the weighted path length $\mathcal{Q}_n$ in a random recursive tree of size $n$ is of order $n^2$. The same is valid for its standard deviation. We conjecture that the sequence $(n^{-2} \mathcal{Q}_n)$ converges in distribution to a non-trivial limit; however, the recursive approach worked out in the proof of Theorem \ref{thm_pathlength}, which also applies to the analysis of the path length in random recursive trees, seems not to be fruitful in this context.

\medskip  \textbf{Outline.} All results are proved in Section \ref{sec:2} starting with the proofs of Theorems \ref{thm_large} and \ref{thm_small} as well as \eqref{Dbound} in Section \ref{sec:typical}. Here, most arguments are based on representations of (weighted) depths as sums of bounded independent random variables which go back to Devroye and Neininger \cite{denefinger04}.
Theorem \ref{thm_binary} and Proposition \ref{propxi} are proved in Section \ref{sec:binary}. In this part, the construction of the limiting process relies on suitable uniform $L_1$-bounds on the increments of the process $\Xi_k(x)_{x \in [0,1]}, k \geq 1,$ while the properties of the limit laws formulated in Proposition \ref{propxi} follow from the distributional fixed-point equation \eqref{fix_process}.
Finally, the proof of Theorem \ref{thm_pathlength} relying on the contraction method is worked out in Section \ref{sec:pathlength}.

\section{Proofs} \label{sec:2}
\subsection{Weighted depths of labelled nodes} \label{sec:typical}

In the permutation model, let $A_{j,k}$ be the event that the node labelled $k$ is in the subtree of the node labelled $j$.
Then, $D_k(n) = \sum_{j=1}^n \I{A_{j,k}} -1 $ and
$W_k(n) = \sum_{j=1}^n j \I{A_{j,k}}$.
It is easy to see that  $A_{1, k}, \ldots, A_{k-1,k}$ and $A_{k+1,k}, \ldots, A_{n,k}$ are two families of independent events; however, there exist subtle dependencies between the sets. Following the approach in \cite{denefinger04}, let $B_{j,k} = A_{j,k-1}$ for $j < k$ and $B_{j,k} = A_{j,k+1}$ for $j > k$. 
For convenience, let $B_{k,k}$ be an almost sure event.
The following lemma summarizes results in \cite{denefinger04} and we refer to this paper for a proof. In this context, 
note that Devroye \cite{devcounter}
gives distributional representations as sums of independent (or $m$-independent) indicator variables for quantities growing linearly in $n$, such as the number of leaves.
\begin{lemma} \label{lem_help}
Let $1 \leq k \leq n$.  Then, the events $B_{j,k}, j = 1, \ldots, n$, are independent. For $j \neq k$, we have
\begin{align*} 
\Prob{A_{j,k}} = \frac{1}{|k-j| + 1}, \quad \Prob{B_{j,k}} = \frac{1}{|k-j|}. 
\end{align*}
\end{lemma}
From the lemma, it follows that
$$\E{ \sum_{j = 1}^n \I{B_{j,k} \backslash A_{j,k}}} \leq 2, \quad \text{and} \quad \E{\sum_{j = 1}^n j \I{B_{j,k} \backslash A_{j,k}}} \leq 2k + \log n.$$
The ideas in \cite{denefinger04} can also be used to analyze second (mixed) moments. Straightforward calculations show the following bounds:
$$\E{ \sum_{i,j = 1}^n \I{B_{j,k}} \I{B_{i,k} \backslash A_{i,k}}} = O(1),$$  and $$\E{ \sum_{i,j = 1}^n i j  \I{B_{j,k}} \I{B_{i,k} \backslash A_{i,k}}} = O(k^2 + k(\log n)^2).$$
Here, both $O$-terms are uniform in $1 \leq k \leq n$.
Define $\bar D_k(n) = \sum_{j=1}^n \I{B_{j,k}}-1$ and $\bar W_k(n) = \sum_{j=1}^n j \I{B_{j,k}}$. We make the following observation:

\begin{itemize} \item [\textbf{O}:]
The asymptotic statements in \eqref{expansionW}, \eqref{expansionV}, Theorem \ref{thm_large} and Theorem \ref{thm_small} are correct if and only if they are correct upon replacing $D_k(n)$ by $\bar D_k(n)$ and/or $W_k(n)$ by $\bar W_k(n)$.
\end{itemize}

 For $i=1,2, n \geq 0$ and $1 \leq k \leq n$, set $H^{(i)}_{n} := \sum_{j=1}^n j^{-i}$ and $H^{(i)}_{k,n} :=  H^{(i)}_{k-1} + H^{(i)}_{n-k}$. Using Lemma \ref{lem_help}, one easily computes 
\begin{align*}  \E{\bar W_k(n)} &=   k (H_{k,n}^{(1)}-1) + n  + 1 ,\\
  \text{Var}(\bar W_k(n)) &=  k^2  (H_{k,n}^{(1)} - H_{k,n}^{(2)}-3)  +  \frac{n^2}{2} + kn + 2k (H^{(1)}_{k-1} - H^{(1)}_{n-k}) - \frac n 2 + k + 1.
\end{align*}
As $H^{(1)}_{n} = \log (n+1) + O(1)$ and $H^{(2)}_n = O(1)$, both expansions \eqref{expansionW} 
and  \eqref{expansionV} follow from  observation \textbf{O}.

\subsubsection{Weighted depths of large nodes}
We prove Theorem \ref{thm_large}. First, \eqref{diff_WD} follows  from \eqref{expansionD} and
\begin{align} \label{bounddk}  \E{\left | k D_k(n) - W_k(n) \right|} \leq k + \sum_{j=1}^n |k-j| \Prob{A_{j,k}} \leq k+n. \end{align}
For $k = \omega(n/\sqrt{\log n})$, combining \eqref{expansionD}, \eqref{limDN} and \eqref{expansionW}, in distribution, 
\begin{align*} \left( \frac{D_k(n) - \E{D_k(n)}}{\sigma_{ D_k(n)}}, \frac{W_k(n) - \E{W_k(n)}}{\sigma_{W_k(n)}} \right) \to (\sN, \sN). \end{align*}
From here, statement \eqref{lim_DW_large} follows from \eqref{expansionD} and \eqref{expansionW}.

Considering the last inserted node with value $Y_n$, note that, conditionally on $Y_n = k$, the correlations between the events $A_{j,k}, j < k$ and $A_{j,k}, j > k$ vanish. More precisely, given $Y_n = k$, the family $\{ \I{A_{j,k}}, j = 1, \ldots, n\}$ is distributed like a family of independent Bernoulli random variables $\{V_{j,k}: j = 1, \ldots, n\}$ with $\Prob{V_{j,k}=1} = |k-j|^{-1}$ for $j \neq k$ and $\Prob{V_{k,k}=1} = 1$. Thus,
\begin{align*}
\E{|Y_n (X_n+1) - \mathbb X_n|} & \leq \frac 1 n \sum_{k=1}^n \E{\sum_{j=1}^n |k-j| \I{A_{j,k}} \Bigg| Y_n = k} \\ 
& = \frac 1 n \sum_{k=1}^n\E{\sum_{j=1}^n |k-j| V_{j,k}} \leq n.
\end{align*}
By \eqref{dev2}, we have $X_n / \log n \to 2$ in probability. Hence, in order to prove \eqref{lim_last}, it suffices to show that, in distribution,
\begin{align} \label{ver} \left(\frac{X_n - 2 \log n}{\sqrt{2 \log n}}, \frac{Y_n}{n} \right) \to \left(\sN, \xi \right). \end{align}
For a sequence $(k_n)$ satisfying $sn \leq k_n \leq tn$ for $0 < s < t < 1$, let us condition on the event $Y_n = k_n$. Then, by the central limit theorem for triangular arrays of row-wise independent uniformly bounded random variables with diverging variance applied to $V_{j, k_n}, j =1, \ldots, n$, in distribution, 
$$\frac{X_n - 2 \log n}{\sqrt{2 \log n}} \to \sN.$$
Hence, \eqref{ver} follows from an application of the theorem of dominated convergence noting that $Y_n$ is uniformly distributed on $\{1, \ldots, n\}$.
\subsubsection{Weighted depths of small nodes}
We prove Theorem \ref{thm_small}. Let $\bar D^>_k(n) = \sum_{j = k+1}^n \I{B_{j,k}}$ and $\bar W^>_k(n) = \sum_{j = k+1}^n j \I{B_{j,k}}$. Since $k = O(n/ \sqrt{\log n})$, the same calculation as in \eqref{bounddk} shows that,   \begin{align} \label{de1} \frac{\E{| \bar W_k(n) - \bar W^>_k(n) - k (\bar D_k(n) - \bar D^>_k(n) )|}}{n} \leq \frac k n \to 0, \quad n \to \infty. \end{align}
For $\lambda, \mu  \in \R$, we have
\begin{align*}
\log & \E{\exp \left(i \lambda (\bar D^>_k(n) - \log n))/ \sqrt{\log n} + i \mu (\bar W^>_k(n) - k \bar D^>_k(n))/n  \right)} \\
& = - i\lambda \sqrt{\log n} + \log \E{\exp \left(i \sum_{j = k+1}^n \left(\frac{\lambda}{\sqrt{\log n}} + \mu \frac{j -k}{n}\right) B_{j,k} \right)}  \\
& = - i\lambda \sqrt{\log n} + \sum_{j=k+1}^n \log \left(1 + \frac{ \exp \left(i \left(\frac{\lambda}{\sqrt{\log n}} +  \mu \frac{j -k}{n} \right)\right)-1}{j-k}  \right).
\end{align*}
By a standard Taylor expansion, the last display equals
\begin{align*}
& -i  \lambda \sqrt{\log n} + \sum_{j=k+1}^n   \frac{ \exp \left(i \left(\frac{\lambda}{\sqrt{\log n}} +  \mu \frac{j -k}{n} \right)\right)-1}{j-k}   + o(1) \\
& = - i\lambda \sqrt{\log n}  + \sum_{j=k+1}^n   \frac{ \exp  \left( i \mu \frac{j -k}{n} \right)\left(1 + \frac{i \lambda}{\sqrt{\log n}} - \frac{\lambda^2}{2\log n}  \right)-1}{j-k}   + o(1) \\
& = - \lambda^2 /2  + \left(1 + \frac{i \lambda}{\sqrt{\log n}} - \frac{\lambda^2}{2 \log n}  \right) \sum_{j=0}^{n-1}  \frac{ \exp  \left( i \mu \frac{j+1}{n} \right)-1}{j+1}   + o(1) \\
& = - \lambda^2/2 +  \int_0^1 \frac{e^{i \mu x}-1}{x} dx + o(1).
\end{align*}
Here, in the last step, we have used that the sum on the right hand side is a Riemann sum over the unit interval whose mesh size $n^{-1}$ tends to zero.
Thus, using the notation of the theorem, \eqref{fourier_dickman} and L{\'e}vy's continuity theorem, in distribution,
\begin{align} \label{e1}\left(\frac{\bar D^>_k(n) - \log n}{\sqrt{\log n}}, \frac{\bar W^>_k(n) - k \bar D^>_k(n)}{n}\right) \to (\sN, \sY). \end{align}
In order to deduce \eqref{lim_small} note that, by Lemma \ref{lem_help},  $\bar D_k(n)  - \bar D^>_k(n)$ and $(\bar D^>_k(n), \bar W^>_k(n))$ are independent while
\begin{align*}  \frac{\bar D_k(n) - \bar D^>_k(n) - \E{\bar D_k(n) - \bar D^>_k(n)}}{\sigma_{\bar D_k(n) - \bar D^>_k(n)}} \to \sN, \end{align*}
in distribution if and only if $k \to \infty$ using the central limit theorem for sums of independent and uniformly bounded random variables. Since
\begin{align*}
\frac{\bar D_k(n) - \E{\bar D_k(n)}}{\sigma_{\bar D_k(n)}}  & =  \frac{\bar D^>_k(n) - \E{D^>_k(n)}}{\sqrt{\log n}} \frac{\sqrt{\log n}}{\sigma_{\bar D_k(n)}}  \\ 
& \; + \frac{\bar D_k(n) - \bar D^>_k(n) - \E{\bar D_k(n) - \bar D^>_k(n)}}{\sigma_{\bar D_k(n) - \bar D^>_k(n)}} \frac{\sigma_{\bar D_k(n) - \bar D^>_k(n)}}{\sigma_{\bar D_k(n)}},  \end{align*}
we deduce
\begin{align*} \left(\frac{\bar D_k(n) - \E{\bar D_k(n)}}{\sigma_{\bar D_k(n)}}, \frac{\bar W^>_k(n) - k \bar D^>_k(n)}{n}\right) \to (\sN, \sY), \end{align*}
from \eqref{e1}
upon treating the cases $k = O(1)$ and $k = \omega(1)$ separately. From here, the assertion \eqref{lim_small} follows with the help of \eqref{de1} and observation \textbf{O}.

\subsubsection{Proof of \eqref{Dbound}}
The main observation is that the $k$th external node visited by the depth first search process is always contained in the subtree rooted at the node labelled $k$.
This can be proved by induction exploiting the decomposition of the tree at the root. Thus, denoting by $H_k(n)$ the height of the subtree rooted at the node labelled $k$, we have
\begin{align*}
D_k(n) & \leq D_k^*(n) \leq D_k(n) + H_k(n), \\
W_k(n) & \leq W_k^*(n) \leq W_k(n) + M_k(n) H_k(n).
\end{align*}
Here, $M_k(n)$ stands for the largest label in the subtree rooted at the node labelled $k$. Let $T_k(n)$ be the size of the subtree rooted at $k$. Then $T_k(n) = 1 + T^{<}_k(n) + T^{>}_k(n)$ where $T^{<}_k(n)$ denotes the number of elements in the subtree rooted at $k$
with values smaller than $k$. By Lemma \ref{lem_help}, for $\ell \leq n - k$, we have $\Prob{T^{>}_k(n) \geq \ell} = \Prob{A_{k,k+\ell}} = 1/(\ell+1)$. Using the same arguments for the quantity $T^{<}_k(n)$, we deduce that,  uniformly in $1 \leq k \leq n$,
\begin{align*}
\E{T_k(n)} = \Theta(\log n), \quad \E{(T_k(n))^2} = \Theta(n^{1/2}), \quad \E{(\log T_k(n))^2} = O(1).
\end{align*}
Thus, by an application of \eqref{conv_height}, for some $C_1 > 0$,
\begin{align*}
\E{|D_k(n) - D_k^*(n)|^2} \leq \E{(H_k(n))^2} & \leq C_1 \E{(\log T_k(n))^2} =  O(1). 
\end{align*}
 By the same arguments,  for some $C_2 > 0$, we have
\begin{align*}
\E{|W_k(n) - W_k^*(n)|^2} & \leq \E{(M_k(n) H_k(n))^2} \leq \E{(k + T_k(n))^2 (H_k(n))^2} \\
& \leq C_2 k^2 + 2 C_1 k \E{T_k(n) (\log T_k(n))^2}  \\ 
 & \, \,  \,  \,+ C_1 \E{(T_k(n))^2 (\log T_k(n))^2}  \\
& = O(k^2 +  (\log n)^{2} n^{1/2}).
\end{align*}
From here, \eqref{Dbound} follows from \eqref{expansionW}.
\subsection{The weighted silhouette} \label{sec:binary} \hfill \\

We prove Theorem \ref{thm_binary} and Proposition \ref{propxi}.

\medskip \textbf{Proof of Theorem \ref{thm_binary}.} We start with the uniform convergence of $(\Xi_k)$. 
For all $x \in [0,1]$, $|\Xi_k(x) - \Xi_{k-1}(x)|$ is distributed like the product of $k+1$ independent random variables, each of which having the uniform distribution on $[0,1]$. In particular, by the union bound and Markov's inequality, for any $m \geq 1$,
$$\Prob{ \sup_{x \in [0,1]} |\Xi_k(x) - \Xi_{k-1}(x)| \geq t} \leq 2^k \Prob{ \prod_{i=1}^{k+1} U_i \geq t} \leq \left(\frac{2}{m+1}\right)^{k} t^{-m}.$$
For $k \geq 1$, let $\mathcal{D}_k = \{\ell 2^{-k}: \ell = 1, \ldots, 2^k - 1\}$.
By construction, for $k \geq 1$, the map $x \to \Xi_k(x)$ is a right-continuous step function. Further, it is continuous at $x$ if and only if $x \notin \mathcal{D}_k$. Next,  for $0 < q < 1$,
\begin{align*}
\E{ \sup_{x \in [0,1]} |\Xi_k(x) - \Xi_{k-1}(x)| } & = \int_0^\infty \Prob{ \sup_{x \in [0,1]} |\Xi_k(x) - \Xi_{k-1}(x)| \geq t} dt \\
& \leq q^k + \int_{q^k}^\infty \left(\frac{2}{m+1}\right)^{k} t^{-m} dt \\
& = q^k +  \frac{1}{m-1}\left(\frac{2}{m+1}\right)^{k} q^{-k(m-1)}.
\end{align*}
With $m=2$ and $q = \sqrt{2/3}$ the latter expression is bounded by $2 q^k$. By Markov's inequality, it follows that $\sup_{m \geq n} \sup_{x \in [0,1]} |\Xi_m(x) - \Xi_{n}(x)| \to 0$ 
in probability as $n \to \infty$.  An application of the triangle inequality shows that
$\sup_{m, p \geq n} \sup_{x \in [0,1]} |\Xi_m(x) - \Xi_{p}(x)| \to 0$ in probability as $n \to \infty$. By monotonicity, this convergence is almost sure.
Thus, almost surely,  $(\Xi_k)$ is uniformly Cauchy in the space 
of c{\`a}dl{\`a}g functions endowed with the uniform topology. 
By completeness,  $(\Xi_k)$ converges to a limit denoted by $\Xi$ with c{\`a}dl{\`a}g paths. Moreover, $\Xi$ is continuous at
$x \notin \mathcal{D}$ where $\mathcal D = \cup_{m \geq 1} \mathcal{D}_m$ since this is true for all $\Xi_k$, $k \geq 1$. For $x \in \mathcal{D}$, let $\Phi(x)$ be the key of the node associated with $x_1 \ldots x_{k-1}$ where $k \geq 1$ is chosen minimal with $x \in \mathcal{D}_k$. Then, 
$\lim_{y \uparrow x} \Xi(x) = \Phi(x) = \Xi(x).$
Thus, $x \mapsto \Xi(x)$ is continuous. By the construction of the tree, it is clear that $\Xi(x) < \Xi(y)$ for any $x, y \in \mathcal{D}$ with 
$x < y$. As $\mathcal{D}$ is dense in $[0,1]$, the process $\Xi$ is strictly monotonically increasing. Obviously, $\Xi(0) = 0$ and $\Xi(1) = 1$; 
hence, $\Xi$ is the distribution function of a probability measure on $[0,1]$. 

We  turn to the convergence of $\mathcal{B}_n(x)$. For any fixed $x \in [0,1]$,
display \eqref{dev2} implies  that, as $n \to \infty$, in probability, $B_n(x) / \log n \to 1$. Thus, \eqref{conv_rep} follows from the convergence $\Xi_k(x) \to \Xi(x)$. 
 The convergence \eqref{conv_rep} is with respect to all moments since $B_n(x) \leq H_n$ and we have convergence of all moments in \eqref{conv_height}. By the theorem of dominated convergence, for any $m \geq 1$,
again using \eqref{conv_height}, we have
$$\int_0^1 \E {\left| \frac{\mathcal B_n(x)}{\log n} - \Xi(x) \right|^m} dx \to 0. $$
This shows \eqref{Lm_conv}. To prove \eqref{weighted_height}, note that, for any $k \geq 1$, $\sup_{x \in [0,1]} \mathcal B_n(x)$ is larger than the product of the height of the subtree rooted at the node $w_k := 1\ldots 1$ on level $k$ 
and $\Xi_{k-1}(\mathbf 1)$. Let $\varepsilon > 0$. Fix $k$ large enough such that $\Prob{\Xi_{k-1}(\mathbf 1) < 1-\varepsilon} < \varepsilon$. Conditional on its size, the subtree rooted at $w_k$ is a random 
binary search tree. Since its size grows linearly in $n$ as $n \to \infty$, it follows from \eqref{conv_height} that, for all $n$ sufficiently large, its height exceeds $(c^*-\varepsilon) \log n$ with probability at least
$1-\varepsilon$. For these values of $n$, we have $\sup_{x \in [0,1]} \mathcal B_n(x) \geq (c^* - 6 \varepsilon) \log n$ with probability at least $1-2 \varepsilon$. 
As $\varepsilon$ was chosen arbitrarily, this shows \eqref{weighted_height}.

For the joint convergence of $B_n(x)$ and $\mathcal{B}_n(x)$ for fixed $x \in [0,1]$,  we abbreviate $B_n := B_n(x), \mathcal{B}_n := \mathcal{B}_n(x)$, $\Xi_k := \Xi_k(x), \Xi = \Xi(x)$ and $\bar B_n = (B_n - \log n)/\sqrt{\log n}$. Note that $\Xi$ and $B_n$ are not independent which causes the proof to be more technical. Denote by $N_k$ the time when the node associated with $x_1 \ldots x_k$ is inserted in the binary search tree. For any $\varepsilon > 0$, we can choose $k, L \geq 1$ such that, for all $n$ sufficiently large,
$$\Prob{|\Xi_k - \Xi| \geq \varepsilon} + \Prob{N_k \geq L} + \Prob{\left | \frac{\mathcal{B}_n}{\log n} - \Xi \right| \geq \varepsilon} \leq \varepsilon.$$
Further, there exists $\delta > 0$ such that $\Prob{|\Xi_k - \Xi_{k-1}| \leq \delta} \leq \varepsilon$.
Then, for $r , y  \in \R$ with $\Prob{\Xi = y } = 0$, and $n$ large enough,
\begin{align*}
\Prob{\bar B_n \leq r, \frac{\mathcal{B}_n}{ \log n} \leq y} \leq 2 \varepsilon + \Prob{\bar B_n \leq r, \Xi_k \leq y + 2 \varepsilon, |\Xi_k - \Xi_{k-1}| \geq \delta, N_k < L}.
\end{align*}
Let $\bar x = x_{k+1} x_{k+2} \ldots$, $(V_1, V_2, \ldots)$ be an independent copy of $(U_1, U_2, \ldots)$ and $$\text{Bin}(n,p) := \sum_{i=1}^n \I{\{V_i \leq p\}}, \quad n \geq 0, p \in [0,1].$$
Given $\Xi_k, |\Xi_k - \Xi_{k-1}|, N_k$, on $N_k < n$, $\bar B_n$ is distributed like $\bar B^*_{\text{Bin}(n - N_k, |\Xi_k - \Xi_{k-1}|)}(\bar x)  + k / \sqrt{\log n}$ where $(B^*_n(\bar x))$ is distributed like $(B_n(\bar x))$ and independent from the remaining quantities. We deduce
\begin{align*}
& \Prob{\bar B_n \leq r, \frac{\mathcal{B}_n}{ \log n} \leq y}  \\ & \leq 2 \varepsilon + \Prob{\frac{k}{\sqrt{  \log n}} + \bar B^*_{\text{Bin}(n - L, \delta)}(\bar x) \leq r, \Xi_k \leq y + 2 \varepsilon, |\Xi_k - \Xi_{k-1}| \geq \delta, N_k < L} \\
& \leq  3 \varepsilon + \Prob{\frac{k}{\sqrt{  \log n}} + \bar B^*_{\text{Bin}(n - L, \delta)}(\bar x) \leq r} \Prob{\Xi \leq y + 2 \varepsilon}
\end{align*}
Using the asymptotic normality of $(\bar B_n^*(\bar x))$ (after rescaling) in \eqref{dev2}, taking the limit superior as $n \to \infty$ and then letting $\varepsilon$ tend to zero, we obtain
$$\limsup_{n \to \infty} \Prob{\bar B_n \leq r, \frac{\mathcal{B}_n}{\log n} \leq y} \leq \Prob{\sN \leq r} \Prob{\Xi \leq y}.$$
The proof of the converse direction establishing \eqref{conv_binary_joint} is easier. It runs along the same lines upon using the trivial bounds $|\Xi_k - \Xi_{k-1}| \leq 1$ and $N_k \geq 0$. 

\medskip \textbf{Proof of Proposition \ref{propxi}.}
We start with the characterization of the distribution of the process. For a deterministic sequence of pairwise different numbers $u_1, u_2, \ldots$ on the unit interval, we define $\xi_k(x)$ analogously to $\Xi_k(x)$ in the infinite binary search tree constructed from this sequence. Here, we abbreviate $\xi_k(x) = 0$ if the node
$x_1 \ldots x_k$ is not in the tree. 
Let $n_m^-, m \geq 1,$ be the subsequence defined by the elements $u_{n^-_m} < u_1$ and $u_m^+, m \geq 1$, be the subsequence defined by the elements $u_{n^+_m} > u_1$. At least one of these sequences is infinite. For $m \geq 1$, let
$y_m^{-} = u_{n^{-}_m} / u_1$ and $y_m^+ = (u_{n^+_m} - u_1) / (1-u_1)$. Next, define $\xi^{-}_k$ ($\xi^+_k$, respectively) analogously to $\xi_k$ based on the sequence $(y^-_m)$ ($(y^+_m)$, respectively). By construction, for $k \geq 1$,  $$\xi_k(x) = \I{[0,1/2)}(x) u_1 \xi_{k-1}^-(2x) + \I{[1/2,1]}(x) ((1-u_1) \xi_{k-1}^+(2x - 1) + u_1).$$
Applying the construction to the sequence $U_1, U_2, \ldots$ yields
$$\Xi_k(x) = \I{[0,1/2)}(x) U_1 \Xi_{k-1}^-(2x) + \I{[1/2,1]}(x) ((1-U_1) \Xi_{k-1}^+(2x - 1) + U_1).$$
Almost surely, the random sequences $y_m^{-}$ and $y_m^+$ are both infinite and $(\Xi^-_k), (\Xi^+_k)$ are independent copies of $(\Xi_k)$. Further, both 
sequences are independent of $U_1$. Hence, letting $k \to \infty$ in the last display, we obtain \eqref{fix_process} on an almost sure level. The characterization of 
$\Law(\Xi)$ by \eqref{fix_process} follows from a standard contraction argument, the argument on page 267 in \cite{grro} applies to our setting without any modifications.

We move on to the statements $i)$ -- $vi)$ on the marginal distributions of the process. Here, we use notation that was introduced in the proof of Theorem \ref{thm_binary}. By continuity, it suffices to show $i)$ 
for $x \in \mathcal{D}$. 
Let $k \geq 1$. By symmetry, for $1 \leq i \leq 2^k-1$, we have $\E{\Phi(i 2^{-k})} = i 2^{-k}$. Thus, the assertion follows for $x \in \mathcal{D}$ since $\Phi(x) = \Xi(x)$. The symmetry statement $ii)$ is reminiscent of the fact that the uniform distribution on $[0,1]$ is symmetric around $1/2$. More precisely, we apply the reflection argument from \cite{aglama07} which is at the core of the proof of the second assertion in \eqref{lasm}. Let $U_1^* = 1 - U_1, U_2^* = 1 - U_2, \ldots$  and define $\Xi^*$ analogously to $\Xi$ in the binary search tree process relying on the sequence $U_1^*, U_2^*, \ldots$ Then, $\Xi^*(t) +\Xi(1-t) = 1$ for all $t \in [0,1]$ which proves $ii)$.
With $Y = \Xi(\xi)$, \eqref{fix_process} yields
$$\Law(Y) = \Law(U Y + \I{A} (1-U)),$$ where $\I{A}, U, Y$ are independent and $\Prob{A} = 1/2$. From \cite{anotherarcsine}, it follows that $Y$ has the arcsine distribution, proving $iii)$.
We move on to the statements about the distribution of $\Xi(t)$. 
Let $t \in (0,1/2)$. Since $\Xi$ is strictly increasing, we have $\Xi(2t) \in (0,1)$ almost surely.  By \eqref{fix_process}, $\Law(\Xi(t)) = \Law(U \Xi(2t))$ with conditions as in \eqref{fix_process}. Therefore, $\Law(\Xi(t))$ admits a density. By symmetry, the same is true for $t \in (1/2,1)$. For $t \in (0,1/2)$, by conditioning on the value of $U$, one finds the density
\begin{align} \label{den1} f_t(x) = \E{ \frac{\I{[x,1]}(\Xi(2t))}{\Xi(2t)}}, \quad x \in (0,1]. \end{align}
 $f_t(x)$ is monotonically decreasing and continuous on $(0,1]$ with $f(1)=0$.  For $t \in (1/2,1)$, $f_t(x) = f_{1-t}(1-x), x \in (0,1)$ is a density of $\Law(\Xi(t))$ by $ii)$. By \eqref{den1}, for $t \in (0,1/2), x \in (0,1)$, 
\begin{align} \label{den2} f_t(x) = \int_x^1 \frac{f_{2t}(y)}{y} dy, \quad \text{or} \quad x f_t'(x) = -f_{2t}(x). \end{align}
Upon setting $f_0 = f_1 = 0$, the last identity also holds for $t =0$ and $t = 1/2$ since $f_{1/2} = \I{[0,1]}$ is a density of $\Law(\Xi(1/2))$. Thus,  for any $t \in (0,1)$, $f_t$ is smooth on $(0,1)$. 
Since the uniform distribution takes values arbitrarily close to one, it follows that, for all $\delta > 0, t \in (0,1)$, we have $\Prob{\Xi(t) > 1 - \delta} > 0$. Hence, for all $t \in (0,1)$, the density $f_t$ is strictly positive on $(0,1)$. 
Thus, for $t \in (0, 1/2)$, $f_t$ is strictly monotonically decreasing.  Summarizing, we have shown $iv)$ and $v)$. For $t \in (0,1/4]$, the assertion $\alpha^{(0)}_t = \infty$ in $vi)$ follows immediately from \eqref{den2} since $\alpha^{(0)}_{2t} > 0$.
Let $1/4 < t < 1/2$. Assume $\alpha^{(0)}_{1-2(1-2t)} < \infty$. Then, $f_{2(1-2t)}(1) < \infty$. By \eqref{den2}, it follows that $f'_{1-2t}(1)$ is finite and hence $f'_{2t}(0)$ is finite. Thus, $f_{2t}(y) / y$ is bounded in a neighhourhood of zero and $\alpha^{(0)}_t < \infty$. For $t > 3/8$, we have $1-2(1-2t) > 1/2$, thus, $\alpha^{(0)}_t < \infty$. Iterating this argument leads to $\alpha_t^{(0)} < \infty$ for all $1/3 < t < 1/2$. 
In order to proceed further, note that, for $t > 1/4$, there exists $k \in \N$, such that, in probability, $\Xi(t) \geq Z := U_1(U_2 + (1-U_2)\prod_{\ell = 1}^kU_{2+\ell}).$
$Z$ admits a density $f_Z$ given by
$$f_Z(x) = 1 + \int_x^1 r(y) dy - x r(x), \quad r(x) = \frac{1}{x^2} \int_0^x \Prob{\prod_{\ell = 1}^kU_{2+\ell} \leq \frac{x-v}{1-v}} dv.$$
Thus, $$\lim_{x \downarrow 0} f_Z(x) = 1 + \int_0^1 r(y) dy < \infty.$$ It follows that $\alpha_t^{(0)} \leq  1 + \int_0^1 r(x) dx < \infty$. Since $\Xi$ is increasing, the function $t \mapsto \alpha_t^{(0)}$ is decreasing. Thus, by monotonicity and continuity, it follows $\alpha^{(0)}_t \uparrow \infty$ as $t \downarrow 1/4$. 
For $t \leq 1/4$, $\alpha^{(0)}_t = \infty$ follows immediately from $\eqref{den2}$ since $\alpha^{(0)}_{2t} < \infty$. For $1/4 < t < 1/2$, the remaining statements about $\alpha_t^{(1)}$ are direct corollaries of the results for $\alpha_t^{(0)}$ since $\alpha^{(1)}_t = \alpha^{(0)}_{1 -2(1-2t)}$. This finishes the proof of $vi)$.

\medskip \textbf{The curvature.} We make a concluding remark about the curvature of $f_t, t \in (0,1/2)$. First, since $x f^{''}_t(x) = - f_{2t}'(x) - f_t'(x)$, for $0 < t \leq 1/4$,  the function $f_t$ is convex. 
From \eqref{den2} it is easy to deduce $f_{1/3}(x) = 2(1-x)$. Since $f_{1/3}'' = f_{1/2}'' = 0$, it is plausible to conjecture that $f_t$ is convex for $t \leq 1/3$ and concave for $1/3 \leq t < 1/2$. Concavity at rational points with small denominator such as $t = 3/8$ or $t = 5/12$ can be verified by hand using \eqref{den2}.

\subsection{Weighted path length and Wiener index} \label{sec:pathlength} 

In order to obtain mean and variance for the weighted path length and the weighted Wiener index, we use  the reflection argument  from the proof of Proposition \ref{propxi} $ii)$. To this end, let $\Ew_n^*$ and $\Ww_n^*$ denote weighted path length and weighted Wiener index in the binary search tree built from the sequence $U_1^* = 1-U_1, U_2^* = 1-U_2, \ldots$ Then, $\Ew_n + \Ew_n^* = P_n + n$ and $\Ww_n + \Ww_n^* =  W_n + {n \choose 2}$ providing the claimed expansions for $\E{\Ew_n}$ and $\E{\Ww_n}$ upon recalling \eqref{meanpath} and \eqref{meanwiener}.

For a finite rooted labelled binary tree $T$, denote by $p(T)$ its path length, by $\ew(T)$ its weighted path length, by $w(T)$ its Wiener index and by $\ww(T)$ its weighted Wiener index. Let $T_1, T_2$ be its left and right subtree and $x$ the label of the root. Then, denoting by $|T|$ the size of $T$, for $|T| \geq 1$,
\begin{align}
p(T) & = p(T_1) + p(T_2) + |T|-1, \label{rec_unlabelled} \\
w(T) & = w(T_1) + w(T_2) + (|T_2| + 1) p(T_1) + (|T_1| + 1) p(T_2) + |T| + 2 |T_1| |T_2| -1.
\end{align}
The first statement is obvious, the argument for the second can be found in \cite{newiener}. For the weighted quantities, one obtains
\begin{align}
\ew(T) & = \ew(T_1) + \ew(T_2) + |T| x,  \\
\ww(T) & = \ww(T_1) + \ww(T_2) + (|T_2| + 1) \ew(T_1) + (|T_1| + 1) \ew(T_2) + (|T| +  |T_1| |T_2|)x. \label{rec_labelled}
\end{align}
Again, the first assertion is easy to see and we only justify the second. The terms $\ww(T_1)$ and $\ww(T_2)$ account for weighted distances within the subtrees. The sum of all weighted distances between nodes in the left subtree and the root equals $\ew(T_1) + |T_1|x$. Replacing $T_1$ by $T_2$, we obtain the analogous sum in the right subtree. The sum of all distances between nodes in different subtrees equals $|T_1| \ew(T_2) + |T_2| \ew(T_1) + |T_1| |T_2| x$. Finally, we need to add $x$ for the weighted distance of the root to itself. Adding up the terms and simplifying leads to \eqref{rec_labelled}.
For $\alpha, \beta > 0$ let $\alpha T + \beta$ be the tree obtained from $T$ where each label $y$ is replaced by $\alpha y + \beta$. Obviously, $p(T) = p(\alpha T +\beta)$ with the analogous identity for the Wiener index. For the weighted quantities, we have
\begin{align}
\ew(\alpha T + \beta ) & = \alpha \ew(T) + (p(T) + |T|) \beta,  \\
\ww(\alpha T + \beta) & = \alpha \ww(T) + (w(T) + |T|(|T| + 1)/2) \beta. \label{conn_factor}
\end{align}
Let $T$ be the  binary search tree of size $n$ in the i.i.d.\ model. Then, given $I_n := \text{rank}(U_1), U := U_1$, in distribution, the trees $\frac{1}{U} T_1$ and $\frac{1}{1-U} T_2 - \frac{U}{1-U}$ are independent binary search trees of size $I_n - 1$ and $n- I_n$, constructed from
independent sequences of uniformly distributed random variables on $[0,1]$. Thus, combining \eqref{rec_unlabelled}--\eqref{conn_factor}, for the vector
$Y_n = (\Ww_n, W_n,\Ew_n,P_n)^T$, we have
\begin{align*}
 Y_n&\stackrel{d}{=} \left[\begin{array}{cccc} U &\quad 0&\quad (n+1-I_n)U&\quad 0 \\0 & 1& 0&
n+1-I_n\\0 & 0& U& 0\\0 & 0& 0& 1
\end{array}\right]Y_{I_n-1}\nonumber\\
& \: \:+ \left[\begin{array}{cccc} 1-U &\quad U&\quad I_n(1-U)&\quad I_nU
\\0 & 1& 0& I_n\\0 & 0& 1-U& U\\0 & 0& 0& 1
\end{array}\right] Y'_{n-I_n} \\
& \: \: +\left(\begin{array}{cccc}
(2n + (n-I_n)(3 I_n + n - 2))U/2\\n-1+2(I_n-1)(n-I_n)\\(2n - I_n)U\\n-1
\end{array}\right),
\end{align*}
where $(Y_n'), (Y_n), (I_n, U)$ are independent and $(Y'_n)$ is distributed like $(Y_n)$. Here, $\stackrel{d}{=}$ indicates that left- and righthand side are identically distributed.

\medskip We consider the sequence $(Z_n)_{n\geq0}$ defined by
\begin{equation*}
Z_n:=\left(\frac{\Ww_n-\mathbb{E}[\Ww_n]}{n^2},\frac{W_n-\mathbb{E}[W_n]}{n^2},\frac{\Ew_n-\mathbb{E}[\Ew_n]}{n},\frac{P_n-\mathbb{E}[P_n]}{n}\right)^{T}, \quad n \geq 1,
\end{equation*}
and $Z_0 = 0$. Let  $\alpha_n = \E{\Ww_n}, \beta_n = \E{W_n} , \gamma_n = \E{ \Ew_n}$ and $\delta_n = \E{P_n}$. Further, let
\begin{align*} A_1^{(n)}&=\left[\begin{array}{cccc} \left(\frac{I_n-1}{n}\right)^2 U &\quad 0&\quad \left(1-\frac{I_n-1}{n}\right)\frac{I_n-1}{n} U&\quad 0 \\0 & \left(\frac{I_n-1}{n}\right)^2& 0&
\left(1-\frac{I_n-1}{n}\right)\frac{I_n-1}{n}\\0 & 0& \frac{I_n-1}{n} U &
0\\0 & 0& 0&\frac{I_n-1}{n}
\end{array}\right],\\
A_2^{(n)}&=\left[\begin{array}{cccc}
\left(1-\frac{I_n}{n}\right)^2(1-U) &\quad
\left(1-\frac{I_n}{n}\right)^2 U &\quad
\frac{I_n}{n}\left(1-\frac{I_n}{n}\right) (1-U)&\quad
\frac{I_n}{n}\left(1-\frac{I_n}{n}\right) U
\\0 & \left(1-\frac{I_n}{n}\right)^2& 0& \frac{I_n}{n}\left(1-\frac{I_n}{n}\right)\\0 & 0& \left(1-\frac{I_n}{n}\right) (1-U)& \left(1-\frac{I_n}{n}\right) U \\0 & 0&
0&1-\frac{I_n}{n}
\end{array}\right], \end{align*}
and $C^{(n)}=(C_1^{(n)}, C_2^{(n)}, C_3^{(n)}, C_4^{(n)})^T$ with
\begin{align*}
  C_1^{(n)} & = \frac{U}{n^2}\alpha_{I_n-1}+\frac{1-U}{n^2}\alpha_{n-I_n}+\frac{U}{n^2}\beta_{n-I_n}+U\frac{(n+1-I_n)}{n^2}\gamma_{I_n-1}\\
& \: \:  +(1-U)\frac{I_n}{n^2}\gamma_{n-I_n} +U\frac{I_n} {n^2}\delta_{n-I_n} + U \frac{2n + (n-I_n)(3 I_n + n - 2)}{2n^2}-\frac{1}{n^2}\alpha_n,\\
C_2^{(n)} & =\frac{1}{n^2}\beta_{I_n-1}+\frac{1}{n^2}\beta_{n-I_n}+\left(1-\frac{I_n-1}{n}\right)\frac{1}{n}\delta_{I_n-1}+\frac{I_n}{n^2}\delta_{n-I_n}\\
& \: \:  +\frac{n-1+2(n-1)(n-I_n)}{n^2}-\frac{1}{n^2}\beta_n,\\
C_3^{(n)} & =\frac{U}{n}\gamma_{I_n-1}+\frac{1-U}{n}\gamma_{n-I_n}+\frac{U}{n}\delta_{n-I_n}+\left( 2 - \frac{I_n}{n} \right) U-\frac{1}{n}\gamma_n,\\
C_4^{(n)} & =\frac{1}{n}\delta_{I_n-1}+\frac{1}{n}\delta_{n-I_n}+ 1 - \frac{1}{n}-\frac{1}{n}\delta_n.
\end{align*}
Then, from the recurrence for $(Y_n)$, it follows
\begin{equation*}
Z_n\stackrel{d}{=}A_1^{(n)}Z_{I_n-1}+A_2^{(n)}Z'_{n-I_n}+C^{(n)}, \quad n \geq 1,
\end{equation*}
where $(Z_n), (Z'_n), (I_n, U)$ are independent and $(Z'_n)$ is distributed like $(Z_n)$. We prove convergence of $Z_n$ in distribution by an application of the contraction method. To this end, note that $I_n/n  \to U$ almost surely by the strong law of large numbers. Thus, with convergence in $L_2$ and almost surely,
\begin{align*} A_1^{(n)}& \to A_1:=\left[\begin{array}{cccc} U^3 &\quad 0&\quad U^2(1-U)&\quad 0 \\0 & U^2& 0&
U(1-U)\\0 & 0& U^2& 0\\0 & 0& 0&U
\end{array}\right],\\
A_2^{(n)}& \to A_2:=\left[\begin{array}{cccc}
(1-U)^3 &\quad U(1-U)^2&\quad U(1-U)^2&\quad U^2(1-U)\\ 0 & (1-U)^2&
0& U(1-U)\\0 & 0& (1-U)^2& U(1-U)\\0 & 0& 0& 1-U
\end{array}\right],\end{align*}
and
 \begin{align*} 
 C^{(n)} \to C :=\left(\begin{array}{cccc}U^2\log{U}+(1-U^2)\log{(1-U)}+U(-14U^2 + 9U + 5)/4 \\2U\log{U}+2(1-U)\log(1-U)+6U(1-U)\\ U^2\ln{U}+(1-U^2)\ln(1-U)+U \\
 2U\ln{U}+2(1-U)\ln{(1-U)} + 1\end{array}\right).
 \end{align*}
For a quadratic matrix $A$, denote by $\|A\|_{\text{op}}$ its spectral radius.
By calculating the eigenvalues of $A_1 A_1^T$ and $A_2 A_2^T$, one checks that $\|A_1\|_{\text{op}} = U$ and $\|A_2\|_{\text{op}} = 1-U$. Thus,
\begin{align*}
\E{\|A_1 A_1^T\|_{\text{op}}}+\E{\|A_2 A_2^T\|_{\text{op}}} \leq \E{\|A_1\|^2_{\text{op}}}+\E{\|A_2\|^2_{\text{op}}} <1.
\end{align*}
Moreover, we have $\Prob{I_n \in \{1, \ldots, \ell  \} \cup \{n\}} \to 0$ for all fixed $\ell$.
Thus, by Theorem 4.1 in \cite{necontrn}, in distribution and with convergence of the first two moments,
we have $Z_n \to (\Ww,W,\Ew,P)$
where $\mathcal{L}(\Ww,W,\Ew,P)$ is the unique
fixed-point of the map:
\begin{align}
T : \mathcal{M}_2^4(0) \longrightarrow\mathcal{M}_2^4(0), \quad T(\mu) = \mathcal{L} \left(A_1 Z +A_2 Z' + C \right), \label{fixT}
\end{align}
with $A_1, A_2, C$ defined above,
 where $Z, Z', U$ are independent and
 $\mathcal{L}(Z)=\mathcal{L}(Z')=\mu$. Here, $\mathcal{M}_2^4(0)$ denotes the set of probability measures on $\R^4$ with finite absolute second moment and zero mean. Variances and covariances can be computed successively using the fixed-point equation, e.g.\ in the following order:
 $\E{P^2}, \E{P W}$, $\E{W^2}, \E{P \Ew},$ $\E{\Ew^2}, \E{\Ew W},$  $\E{P \Ww}, \E{W \Ww}$, $\E{\Ew \Ww}, \E{\Ww^2}$. Additionally to the variances given in the theorem, one obtains
 \begin{align}
 \text{Cov}(P_n, \Ew_n) & \sim \frac{21 - 2 \pi^2}{ 6} n^2, \quad \text{Cov}(P_n, W_n) \sim \frac{20 -2 \pi^2}{3} n^3, \label{covs}\\
 \text{Cov}(\Ew_n, W_n) & \sim   \frac{10 -  \phantom{2} \pi^2}{3} n^3, \quad \text{Cov}(P_n, \Ww_n) \sim  \frac{10 -  \phantom{2} \pi^2}{3} n^3, \\
 \text{Cov}(W_n, \Ww_n) & \sim  \frac{10 - \phantom{2}  \pi^2}{3} n^4, \quad  \text{Cov}(\Ew_n, \Ww_n)  \sim \frac{481 -48 \pi^2}{288} n^3. \label{cove}
 \end{align}

\section*{Acknowledgements}
The first author is  grateful to the King Saud University, Deanship
of Scientific Research, College of Science Research
Center. The research of the third author was supported by a Feodor Lynen Fellowship of the Alexander von Humboldt-Foundation.

\end{document}